\documentclass[preprint,12pt]{elsarticle}

\usepackage{fancyhdr}
\usepackage{fullpage}
\usepackage{amsmath}
\usepackage{amsbsy}
\usepackage{amssymb}
\usepackage{amscd}
\usepackage{amsfonts}
\usepackage{supertabular}
\usepackage{graphics}
\usepackage{verbatim}
\usepackage{epsfig}
\usepackage{xspace}
\usepackage{euscript}
\usepackage{alltt}
\usepackage{boxedminipage}
\usepackage{float}
\usepackage[colorlinks]{hyperref}
\usepackage{color}
\usepackage[all]{xy}
\usepackage{t1enc}
\usepackage{times,exscale}
\usepackage{graphicx,calc}
\usepackage{subfig}
\usepackage[ruled,vlined]{algorithm2e}
\usepackage{epstopdf}
\usepackage{array,multirow}

\def\br{{\mathbf{r}}}

\def\bF{{\mathbf{F}}}
\def\bA{{\mathbf{A}}}

\def\bX{{\mathbf{X}}}

\def\bff{{\mathbf{f}}}

\newcommand{\norm}[1]{\| #1 \|}

\def\be{{\mathbf{e}}}
\def\br{{\mathbf{r}}}

\def\bx{{\mathbf{x}}}
\def\bb{{\mathbf{b}}}
\def\bA{{\mathbf{A}}}
\def\bI{{\mathbf{I}}}

\journal{arXiv}

\begin{document}

\begin{frontmatter}

\title{Alternating Anderson-Richardson method: An efficient alternative to preconditioned Krylov methods for large, sparse linear systems}
\author[gatech]{Phanish Suryanarayana\corref{cor}}
\author[gatech]{Phanisri P. Pratapa}
\author[llnl]{John E. Pask}
\address[gatech]{College of Engineering, Georgia Institute of Technology, Atlanta, GA 30332, USA}
\address[llnl]{Physics Division, Lawrence Livermore National Laboratory, Livermore, CA 94550, USA}
\cortext[cor]{Corresponding Author (\it phanish.suryanarayana@ce.gatech.edu) }

\begin{abstract}
We present the Alternating Anderson-Richardson (AAR) method: an efficient and scalable alternative to preconditioned Krylov solvers for the solution of large, sparse linear systems on high performance computing platforms. Specifically, we generalize the recently proposed Alternating Anderson-Jacobi (AAJ) method (Pratapa et al., \textit{J.~Comput.~Phys.}~(2016),~306,~43--54) to include preconditioning, discuss efficient parallel implementation, and provide serial MATLAB and parallel C/C++ implementations. 
In serial applications to nonsymmetric systems, we find that AAR is comparably robust to GMRES, using the same preconditioning, while often outperforming it in time to solution; and find AAR to be more robust than Bi-CGSTAB for the problems considered. 
In parallel applications to the Helmholtz and Poisson equations, we find that AAR shows superior strong and weak scaling to GMRES, Bi-CGSTAB, and Conjugate Gradient (CG) methods, using the same preconditioning, with consistently shorter times to solution at larger processor counts. Finally, in massively parallel applications to the Poisson equation, on up to 110,592 processors, we find that AAR shows superior strong and weak scaling to CG, with shorter minimum time to solution. We thus find that AAR offers a robust and efficient alternative to current state-of-the-art solvers, with increasing advantages as the number of processors grows.
\end{abstract}


\end{frontmatter}




\section{Introduction}
Linear systems of equations are encountered in the gamut of applications areas within computational physics, from quantum to continuum to celestial mechanics. The strategies adopted for solving such systems can be broadly classified into two categories: direct methods \cite{davis2006direct} and iterative methods \cite{Saad_text}. For relatively small system sizes, direct methods  such as QR decomposition and LU factorization are generally the preferred approaches. However, as the size of the system increases, direct methods become inefficient---in terms of both computational cost and storage requirements---due to poor scaling with system size relative to iterative approaches, particularly Krylov subspace based techniques such as the Generalized Minimal Residual (GMRES) \cite{saad1986gmres} and Conjugate Gradient (CG) \cite{shewchuk1994introduction} methods. Therefore, such iterative approaches are often the methods of choice for the solution of large-scale linear systems of equations. 

A number of physical applications require the repeated solution of large, sparse linear systems. For example, in real-space quantum molecular dynamics calculations  \cite{jing1994ab,shimojo2005embedded,osei2014accurate,suryanarayana2017SQDFT} or electronic structure calculations with exact exchange \cite{perdew1996rationale,lin2016adaptively}, the Poisson equation may be solved hundreds of thousands of times within a single simulation. Therefore, it is critical to reduce the time to solution as far as possible in such situations, a goal typically achieved through parallel computing, wherein the number of floating point operations per second increases linearly with the number of processors. However, the cost associated with inter-processor communication, especially global communication, limits the parallel efficiency of linear solvers, which in turn limits the reduction in wall time that can be achieved in practice \cite{de1994communication,duff1999developments,yang2003improved}. Therefore, there is wide interest in developing algorithms that scale well on modern large-scale parallel computers which can contain tens to hundreds of thousands of computational cores or more \cite{zuo2010improved,ghysels2013hiding,mcinnes2014hierarchical}.

Krylov subspace methods such as GMRES and CG have limited parallel scalability due to the large number of global operations inherent to them \cite{de1995reducing,ghysels2014hiding}. In this context, the classical Richardson and Jacobi fixed-point iterations \cite{hackbusch1994iterative,Saad_text} are ideally suited for massive parallelization by virtue of the strict locality of operations required, i.e., they do not require the calculation of any dot products, other than those required for the calculation of the residual \cite{golub2001closer,barrett1994templates}. However, they suffer from extremely large prefactors and poor scaling with system size compared to Krylov subspace methods, which has made them unattractive on even the largest modern platforms. This has motivated the development of strategies that significantly accelerate the convergence of the basic Richardson/Jacobi iterations while maintaining their underlying parallel scalability and simplicity to the maximum extent possible \cite{Mittal2014,pratapa2016anderson}. Such approaches include the Chebyshev acceleration technique \cite{Saad_text} and the recently developed Scheduled Relaxation Jacobi (SRJ) method \cite{Mittal2014}. However, Chebyshev acceleration requires the computation of the extremal eigenvalues of the coefficient matrix, which may be computationally expensive. Moreover, the current formulation of the SRJ method is applicable only to linear systems arising from the discretization of elliptic equations using second-order finite-differences. For such reasons, Krylov subspace methods have remained the methods of choice in general for the solution of large, sparse linear systems.

Recently, we proposed to employ Anderson extrapolation \cite{Anderson1965}\footnote{
Anderson's extrapolation has been successfully utilized for accelerating the convergence of non-linear fixed-point iterations arising in electronic structure calculations \cite{pulay1980convergence}, coupled fluid-structure transient thermal problems \cite{ganine2013nonlinear}, as well as neutronics and plasma physics \cite{willert2014leveraging}. In the context of linear systems of equations, Anderson's technique bears a close connection to the GMRES method \cite{rohwedder2011analysis, walker2011anderson, potra2013characterization}.} 
at periodic intervals within the classical Jacobi iteration, resulting in the so called Alternating Anderson-Jacobi (AAJ) method \cite{pratapa2016anderson}. This strategy was found to accelerate the Jacobi iteration by orders of magnitude, to the point in fact of outperforming GMRES significantly in serial computations without preconditioning\footnote{In the context of electronic structure calculations, the analogue of the AAJ method for nonlinear fixed-point iterations---referred to as Periodic Pulay \cite{banerjee2016periodic}---is found to significantly accelerate the convergence of the self-consistent field (SCF) method.}. 
In the present work, we generalize the AAJ method to include preconditioning, discuss efficient parallel implementation, and provide serial MATLAB and parallel C/C++ implementations. 
In serial applications to nonsymmetric systems, we find that AAR is comparably robust to GMRES, using the same preconditioning, while often outperforming it in time to solution.
In parallel applications to the Helmholtz and Poisson equations, on up to 110,592 processors, we find that AAR shows superior strong and weak scaling to GMRES, Bi-CGSTAB, and Conjugate Gradient (CG) methods, using the same preconditioning, with consistently shorter times to solution at larger processor counts. 
We thus find that AAR offers a robust and efficient alternative to current state-of-the-art solvers, with increasing advantages as the number of processors grows.

The remainder of this paper is organized as follows. In Section~\ref{Sec:AAR}, we describe the preconditioned AAR method. We demonstrate the efficiency and parallel scaling of the method in Section~\ref{Sec:Examples}. Finally, we conclude in Section~\ref{Sec:Conclusions}. 


\section{Alternating Anderson-Richardson method} \label{Sec:AAR}
In this section, we present the preconditioned Alternating Anderson-Richardson (AAR) method for the solution of large, sparse linear systems: 
\begin{eqnarray} 
& & \bA \bx = \bb \,, \label{Eqn:LinearSystem} \\
& & \bA \in \mathbb{C}^{N\times N}\,, \,\,\, \bx \in \mathbb{C}^{N\times 1}\,, \,\,\, \bb \in \mathbb{C}^{N\times 1}\,, \nonumber
\end{eqnarray}
where $\mathbb{C}$ is the set of complex numbers. This approach generalizes the Alternating Anderson-Jacobi (AAJ) method presented previously \cite{pratapa2016anderson} to include preconditioning and therefore accelerate convergence. In this work, we summarize and focus on the incorporation of preconditioning and the  development of an efficient parallel formulation and implementation; a more complete discussion of the underlying alternating Anderson approach found in our previous work \cite{pratapa2016anderson}.

In the AAR method, Anderson extrapolation \cite{Anderson1965} is performed periodically within a preconditioned Richardson fixed-point iteration \cite{Saad_text} to accelerate its convergence, while maintaining its parallel scalability to the maximum extent possible. Since the method makes no assumptions about the symmetry of $\bA$, it is applicable to symmetric and nonsymmetric systems alike. Moreover, it is amenable to the three types of preconditioning: left, right, and split \cite{Saad_text,benzi2002preconditioning}. For the sake of simplicity, we choose left preconditioning in the present work. Mathematically, the linear system in Eq.~\ref{Eqn:LinearSystem} is solved using a fixed-point iteration of the form 
\begin{equation} \label{Eqn:UpdateAAR}
\bx_{k+1} = \bx_k + \mathbf{B}_k \bff_k \,, \quad k=0, 1, \ldots
\end{equation}
where the matrix $\mathbf{B}_k \in \mathbb{C}^{N \times N}$ can be written as
\begin{equation} \label{Eqn:Bk}
\mathbf{B}_k = \begin{cases} 
                   \omega \bI & \mbox{if } (k+1)/p \not\in \mathbb{N} \,,\quad\text{(Richardson)}\\
                   \beta \mathbf{I} - (\mathbf{X}_k + \beta \mathbf{F}_k)(\mathbf{F}_k^{T}\mathbf{F}_k)^{-1} \mathbf{F}_k^{T}& \mbox{if } (k+1)/p \in \mathbb{N} \,.\quad\text{(Anderson)}
                   \end{cases}
\end{equation}
Above, $\omega \in \mathbb{C}$ is the relaxation parameter used in the Richardson update, $\beta \in \mathbb{C}$ is the relaxation parameter used in the Anderson update, $\bI \in \mathbb{R}^{N\times N}$ is the identity matrix, the superscript $T$ denotes the conjugate transpose, and $p$ is the frequency of Anderson extrapolation. Additionally, $\bX_k \in \mathbb{C}^{N\times m}$ and $\bF_k \in \mathbb{C}^{N\times m}$ contain the iteration and residual histories at the $k^{th}$ iteration:
\begin{eqnarray}
\bX_k & = & \begin{bmatrix} (\bx_{k-m+1} - \bx_{k-m}) & (\bx_{k-m+2} - \bx_{k-m+1}) & \ldots & (\bx_k-\bx_{k-1}) \end{bmatrix}  \,, \label{Eq:Xk}\\
\bF_k & = & \begin{bmatrix} (\bff_{k-m+1} - \bff_{k-m}) & (\bff_{k-m+2} - \bff_{k-m+1}) & \ldots & (\bff_k-\bff_{k-1}) \end{bmatrix} \,, \label{Eq:Fk}
\end{eqnarray}
where $m+1$ is the number of iterates used for Anderson extrapolation,\footnote{In the initial iterations, $\bx_j$ in Eq.~\ref{Eq:Xk} and $\bff_j$  in  Eq.~\ref{Eq:Fk} with $j<0$ are omitted or can be set to null vectors if a pseudoinverse is used to evaluate $(\mathbf{F}_k^{T}\mathbf{F}_k)^{-1}$ in Eq.~\ref{Eqn:Bk}.} and the residual vector
\begin{equation} \label{Eqn:Residual}
\bff_k = \mathbf{M}^{-1} (\bb - \bA \bx_k) \,,
\end{equation}
with preconditioner $\mathbf{M}^{-1}\in \mathbb{C}^{N\times N}$. As discussed in the context of AAJ \cite{pratapa2016anderson}, the key to the robustness and efficiency of the method is the Anderson extrapolation step which minimizes the $l_2$ norm of the residual in the column space of $\bF_k$, yielding consistent and substantial reductions with increasing history length $m$;\footnote{In practice, $m$ is limited by available memory and/or finite precision effects as the matrix $\mathbf{F}_k^{T}\mathbf{F}_k$ in Eq.~\ref{Eqn:Bk} becomes ill-conditioned.} while the key to the parallel scalability of the method is that the extrapolation is performed only periodically, thus reducing nonlocal communications significantly.

\begin{figure}[h!]
\centering
\includegraphics[keepaspectratio=true,width=0.6\textwidth]{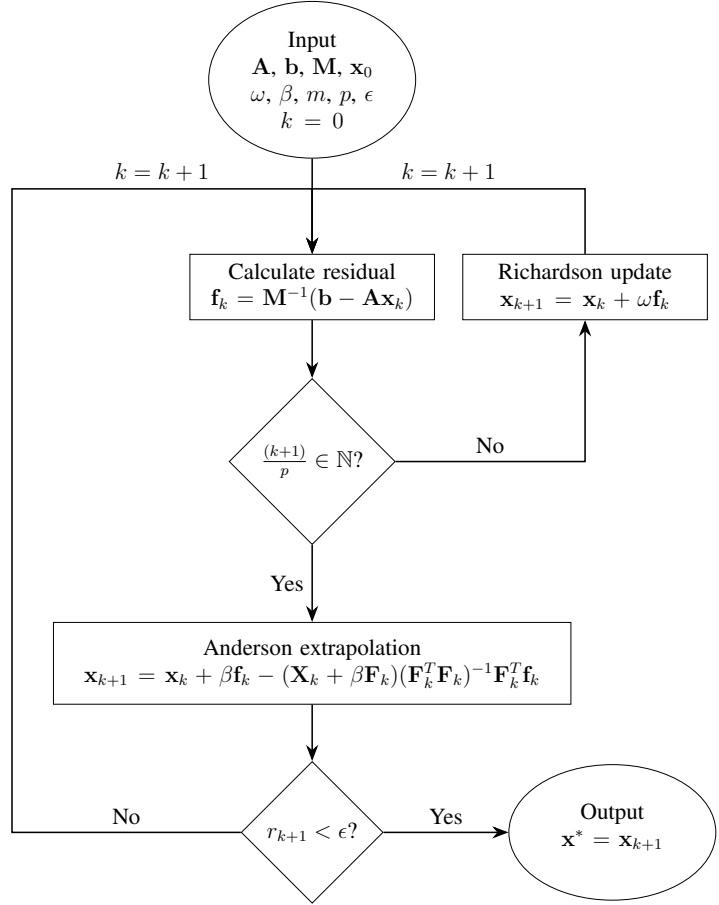} 
\caption{The preconditioned Alternating Anderson-Richardson (AAR) method.}
\label{Fig:flowchart}
\end{figure}

We summarize the AAR method in Fig.~\ref{Fig:flowchart}, wherein $\bx_0$ denotes the initial guess, $r_k$ denotes the relative $l_2$ norm of the residual vector (i.e., $r_{k} = \norm{\bA \bx_k - \bb}/\norm{\bb}$), and $\epsilon$ is the tolerance specified for convergence. In order to enhance parallel scalability, we reduce global communications by checking for convergence of the fixed-point iteration (i.e., calculating $r_k$) only during Anderson extrapolation steps. 

The key difference between the AAR and AAJ \cite{pratapa2016anderson} methods lies in the choice of residual vector Eq.~\ref{Eqn:Residual}. In the AAR method, any available preconditioner $\mathbf{M}^{-1}$ can be employed, whereas in AAJ, $\mathbf{M} = \mathbf{D}$ is the diagonal part of $\mathbf{A}$. The AAR method thus generalizes the AAJ method in the sense that the AAJ method is recovered for the particular choice of preconditioner $\mathbf{M}^{-1} = \mathbf{D}^{-1}$, i.e., the classical Jacobi preconditioner. 
Furthermore, just as the AAJ method can be understood as a generalization of the Jacobi \cite{Saad_text} and Anderson-Jacobi (AJ) \cite{Anderson1965,pratapa2016anderson} methods, the AAR method can be understood as a generalization of the Richardson \cite{Saad_text} and Anderson-Richardson (AR) \cite{rohwedder2011analysis, walker2011anderson, potra2013characterization} methods. Specifically, the AR method is recovered for $p=1$, while the Richardson iteration is obtained in the limit $p \rightarrow \infty$.

As discussed in the context of AAJ \cite{pratapa2016anderson}, the convergence of the AAR method can be understood though its connection to GMRES. First, we note that with complete history (i.e., $m=\infty$), AAR is equivalent to AR for $\omega \neq 0$ and $p \geq 1$ since, upon extrapolation, the residual norm is minimized over the same Krylov subspace regardless of previous extrapolations (in exact arithmetic).\footnote{Excluding potential differences in stagnation \cite{pasini2018}.}
Second, it has been shown \cite{rohwedder2011analysis, walker2011anderson, potra2013characterization} that AR with complete history is equivalent to GMRES without restart, in the sense that the iterates of one can be readily obtained from those of the other (in exact arithmetic).\footnote{Excluding potential differences in stagnation \cite{potra2013characterization}.}
Hence, in the above sense, AAR with complete history is equivalent to GMRES without restart and so must show corresponding convergence.

With finite history and restarts, however, as typical in practice to reduce storage and/or orthogonalization costs, the convergence rates of both AAR and GMRES are generally reduced. And in this context, as demonstrated in Section~\ref{Sec:Examples}, we typically find shorter times to solution for AAR than for GMRES, with increasing advantages for AAR in parallel calculations as the number of processors grows. As discussed in the context of AAJ \cite{pratapa2016anderson}, this may be due in part to the fact that AAR retains and minimizes over the most recent $m$-vector history at each extrapolation, while GMRES begins anew at each restart. The key advantage of AAR, however, in parallel calculations in particular, is that the majority of iterations are simple, computationally local Richardson iterations, with Anderson extrapolations only every $p$ iterations. A study of the mathematical properties of AAR in relation to GMRES and AR can be found in the recent work of Lupo Pasini \cite{pasini2018}.

Finally, in finite precision, other considerations come into play. For example, while for complete history and exact arithmetic, the iterates produced by AAR upon extrapolation are independent of $\omega$ and $p$, this no longer holds with finite history and floating point arithmetic. Nevertheless, as shown in the context of AAJ \cite{pratapa2016anderson}, the dependence is generally weak over a broad range of values so that the method is generally insensitive to the particular choice of values within the range. Similar insensitivity is found for the Anderson extrapolation parameter $\beta$, though larger values can accelerate convergence in better-conditioned (or well preconditioned) problems, consistent with findings in the nonlinear context \cite{banerjee2016periodic}. Finally, while in exact arithmetic, a larger history length $m$ must generally improve convergence, by providing a larger subspace over which to minimize, in finite precision, the increasing condition number of the matrix $\mathbf{F}_k^{T}\mathbf{F}_k$ in Eq.~\ref{Eqn:Bk} with increasing $m$ limits the effective history length in practice. 
Given the general insensitivity of the method to the particular choice of parameter values, we use the same default set $\{\omega,\beta,m,p \} = \{0.6, 0.6, 9, 8\}$ for all systems in the present work. While possible to optimize for a particular application area, we have found these to be sufficient in a broad range of applications, with available preconditioning in particular, as demonstrated in Section~\ref{Sec:Examples}.


\section{Results and discussion} \label{Sec:Examples}
In this section, we demonstrate the efficiency and scaling of the preconditioned Alternating Anderson-Richardson (AAR) method in the solution of large, sparse linear systems of equations. Specifically, we consider an assortment of nonsymmetric systems from Matrix Market\footnote{\url{http://math.nist.gov/MatrixMarket/}} as well as Poisson and complex-valued Helmholtz equations arising in real-space electronic structure calculations, and use the default parameters $\{\omega,\beta,m,p\} = \{0.6, 0.6, 9, 8\}$ in AAR for all systems.
\subsection{Matrix Market: assortment of nonsymmetric systems}
In order to demonstrate the robustness and efficiency of AAR, we first study the relative performance of AAR, GMRES, and Bi-CGSTAB in MATLAB\footnote{\url{https://www.mathworks.com/}} for nonsymmetric linear systems available in the Matrix Market repository.\footnote{The MATLAB implementation of AAR is available as part of the code accompanying this paper. For GMRES and Bi-CGSTAB, the inbuilt functions in MATLAB are utilized.} Specifically, we consider ten matrices that arise in various areas of computational physics, including oil reservoir modeling, fluid dynamics, and the study of plasmas. We use the default MATLAB parameters for GMRES and Bi-CGSTAB, with a vector of all ones as the starting guess $\bx_0$ in all cases. The simulations are performed on a workstation with the following configuration: Intel Xeon Processor E3-1220 v3 (Quad Core, 3.10GHz Turbo, 8MB), 16GB (2x8GB) 1600MHz DDR3 ECC UDIMM. 

In Table~\ref{tab:matmark}, we present the timings (in seconds) obtained for achieving a convergence tolerance of $\epsilon = 10^{-6}$ for two cases: (i) Jacobi preconditioner and (ii) ILU(0) preconditioner. 
We first note that the simple Jacobi preconditioner is insufficient to obtain convergence for a number of these  systems, though AAR is able to converge more systems than GMRES and Bi-CGSTAB. 
On the other hand, we see that ILU(0) preconditioning is sufficient to obtain convergence for all ten systems for AAR, all but one for GMRES, but only half the systems for Bi-CGSTAB.
Moreover, we see that when ILU(0) preconditioning is employed, AAR outperforms GMRES significantly for all but one system; while Bi-CGSTAB can outperform both AAR and GMRES significantly when it converges, but even then outperforms AAR in only two of those five cases.
We thus find that AAR shows comparable robustness to GMRES while often outperforming it, while Bi-CGSTAB, though highly competitive when it converges, shows considerably less robustness in the applications considered.

Importantly, we note that while the default parameters for AAR work well for a wide range of applications, as demonstrated above, the flexibility to tune the parameters can be leveraged to tailor AAR for particular applications areas. For example, by simply tuning $\beta$, the solution times for the utm3060, utm1700b, and sherman5 applications above can be brought down to 0.087 s, 0.026 s, and 0.010 s for $\beta$ = 1.7, 1.3, and 0.8, respectively. More importantly, however, by virtue of the superior parallel scaling of AAR, times to solution can be brought down substantially further still relative to standard solvers such as GMRES and Bi-CGSTAB, as we now show.

\begin{table}[h]
\centering
\begin{tabular}{| c |c | c c c | c c c |}
\hline
Matrix & $N$ & \multicolumn{3}{c|}{Jacobi preconditioner} & \multicolumn{3}{c|}{ILU(0) preconditioner} \\
& & AAR & GMRES & Bi-CGSTAB & AAR & GMRES & Bi-CGSTAB \\
\hline
utm3060 & 3060 & 15.104 & - & 1.805 & 0.283 & 0.191 & 0.043  \\
utm1700a & 1700 &  - & - & - & 0.004 & 0.013 & 0.006 \\
utm1700b & 1700 & 4.733 & - & - & 0.045 & 0.159 & 0.023 \\
fidap029 & 2870 & 0.004 & 0.019 & 0.007 & 0.005 & 0.016 & 0.008 \\
sherman5 & 3312 & 0.028 & 0.095 & 0.016 & 0.014 & 0.025 & 0.014 \\
mcfe & 765 & - & 0.392 & - & 0.007 & 0.018 & -  \\
memplus & 17758 & 0.521 & 1.344 & - & 0.735 & 2.185 & -  \\
add32 & 4960 & 0.021 & 0.070 & - & 0.021 & 0.061 & -  \\
mcca & 180 & 0.019 & 0.023 & 0.008 & 0.013 & 0.021 & -  \\
fs\_680\_3 & 680 & 0.211 & - & - & 0.004 & - & -  \\
\hline
\end{tabular}
\caption{Time taken in seconds by AAR, GMRES, and Bi-CGSTAB in MATLAB for nonsymmetric linear systems from Matrix Market. The symbol `-' is used to indicate that convergence was not achieved within $1000$ sec.}
\label{tab:matmark}
\end{table}


\subsection{Orbital-free Density Functional Theory: Helmholtz equation} \label{sec:helmholtz}
Next, we study the relative performance of AAR, GMRES with standard restarts \cite{saad1986gmres}, GMRES with augmented restarts \cite{baker2005technique} (LGMRES), and Bi-CGSTAB \cite{van1992bi} in PETSc \cite{Petsc1,Petsc2}\footnote{The PETSc implementation of AAR for complex-valued systems is available as part of the code accompanying this paper. For GMRES, LGMRES, and Bi-CGSTAB, the inbuilt functions in PETSc are utilized.}. We consider the periodic Helmholtz problem arising in real-space orbital-free Density Functional Theory (OF-DFT) calculations \cite{Choly2002,ghosh2016higher,Suryanarayana2014524}:
\begin{equation} 
-\frac{1}{4\pi} \nabla^2 V(\br)  + Q \, V(\br) = P \, \rho^{\alpha}(\br) \,\,\, \text{in} \,\,\, \Omega,\quad \begin{cases} V(\br)=V(\br+L_i\hat{\be}_i) \,\,\, \text{on} \,\,\, \partial\Omega \,, \\ \hat{\be}_i\cdot \nabla V(\br)=\hat{\be}_i\cdot \nabla V(\br+L_i\hat{\be}_i) \,\,\, \text{on} \,\,\, \partial\Omega \,, \end{cases}  \label{Eqn:Helmholtz_p} 
\end{equation}
where $V(\br)$ is the kernel potential \cite{wang1998orbital,wang1999orbital}, $\rho(\br)$ is the electron density, $\alpha=\frac{5}{6}+\frac{\sqrt{5}}{6}$, $P=0.003277-i0.009081$, $Q=-0.134992-i0.070225$, $i=\sqrt{-1}$, and $\Omega$ is a cuboidal domain with side lengths $L_i$, unit vectors $\hat{\be}_i$ along each edge, and boundary $\partial \Omega$. The equation is discretized using sixth-order accurate finite-differences on a uniform grid with mesh-size $h$. The computations are parallelized by decomposing the domain into cubical subdomains of equal size, with communication between processors handled via Message Passing Interface (MPI) in the PETSc framework. 

In the Anderson extrapolation step of AAR, we perform only one global communication call (i.e., \texttt{MPI\_Allreduce}) to simultaneously determine $r$ and the complete matrix $\mathbf{F}_k^{T}\mathbf{F}_k$. Since the matrix $\mathbf{F}_k^{T}\mathbf{F}_k$ is generally ill-conditioned, we compute its inverse using the Moore-Penrose pseudoinverse \cite{laub2005matrix}. We employ the default PETSc parameters for GMRES, LGMRES, and Bi-CGSTAB. In all the simulations, we use a vector of all zeros as the starting guess $\bx_0$ and a convergence tolerance of $\epsilon = 10^{-6}$ on the relative residual. We perform the calculations on a computer cluster consisting of $16$ nodes with the following configuration: Altus 1804i Server - 4P Interlagos Node, Quad AMD Opteron 6276, 16C, 2.3 GHz, 128GB, DDR3-1333 ECC, 80GB SSD, MLC, 2.5" HCA, Mellanox ConnectX 2, 1-port QSFP, QDR, memfree, CentOS, Version 5, and connected through InfiniBand cable. 

\begin{figure}[h]
\centering
\subfloat[Strong scaling]{\label{Fig:strong:h}\includegraphics[width=0.47\textwidth]{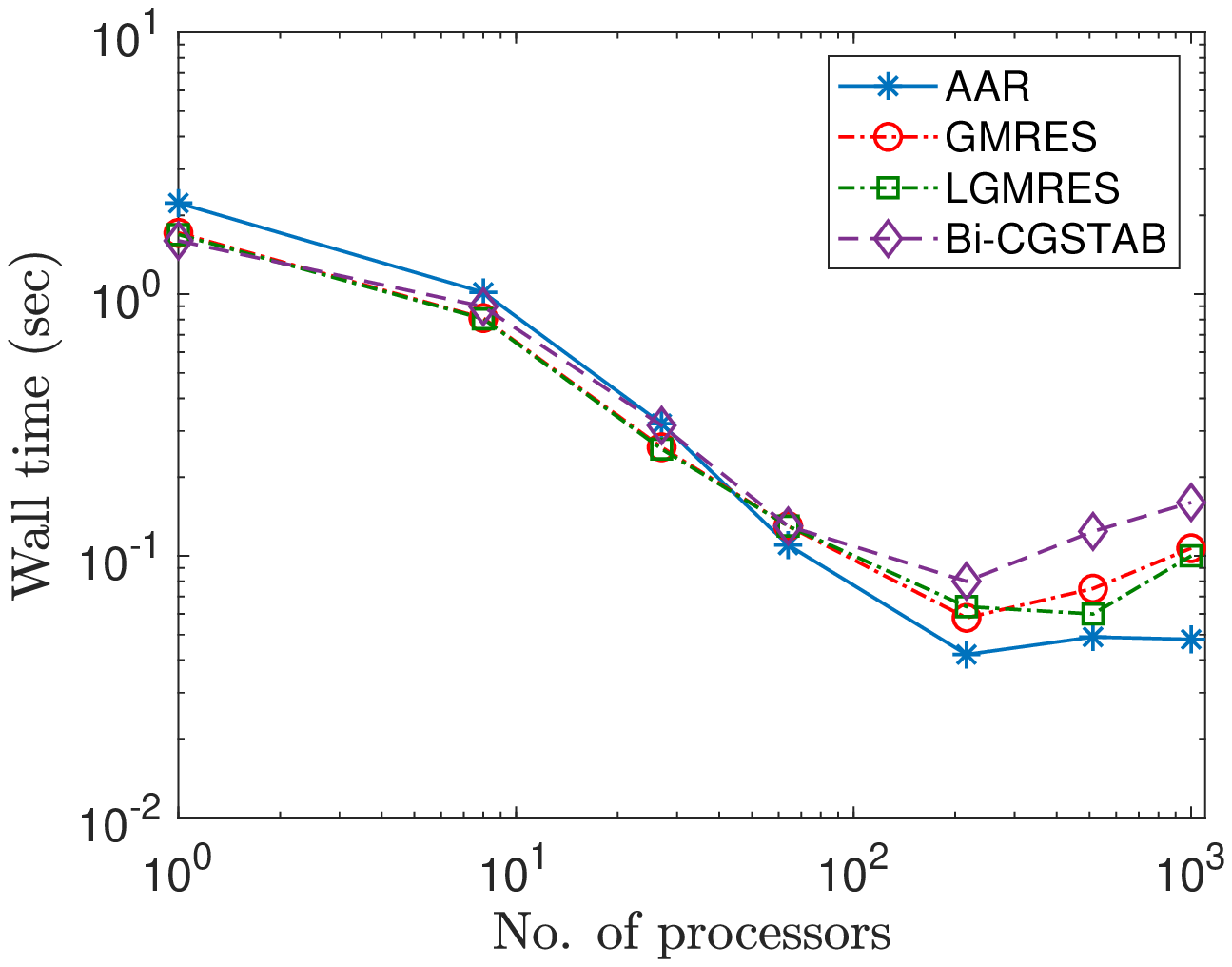}} 
\subfloat[Weak scaling]{\label{Fig:weak:h}\includegraphics[width=0.47\textwidth]{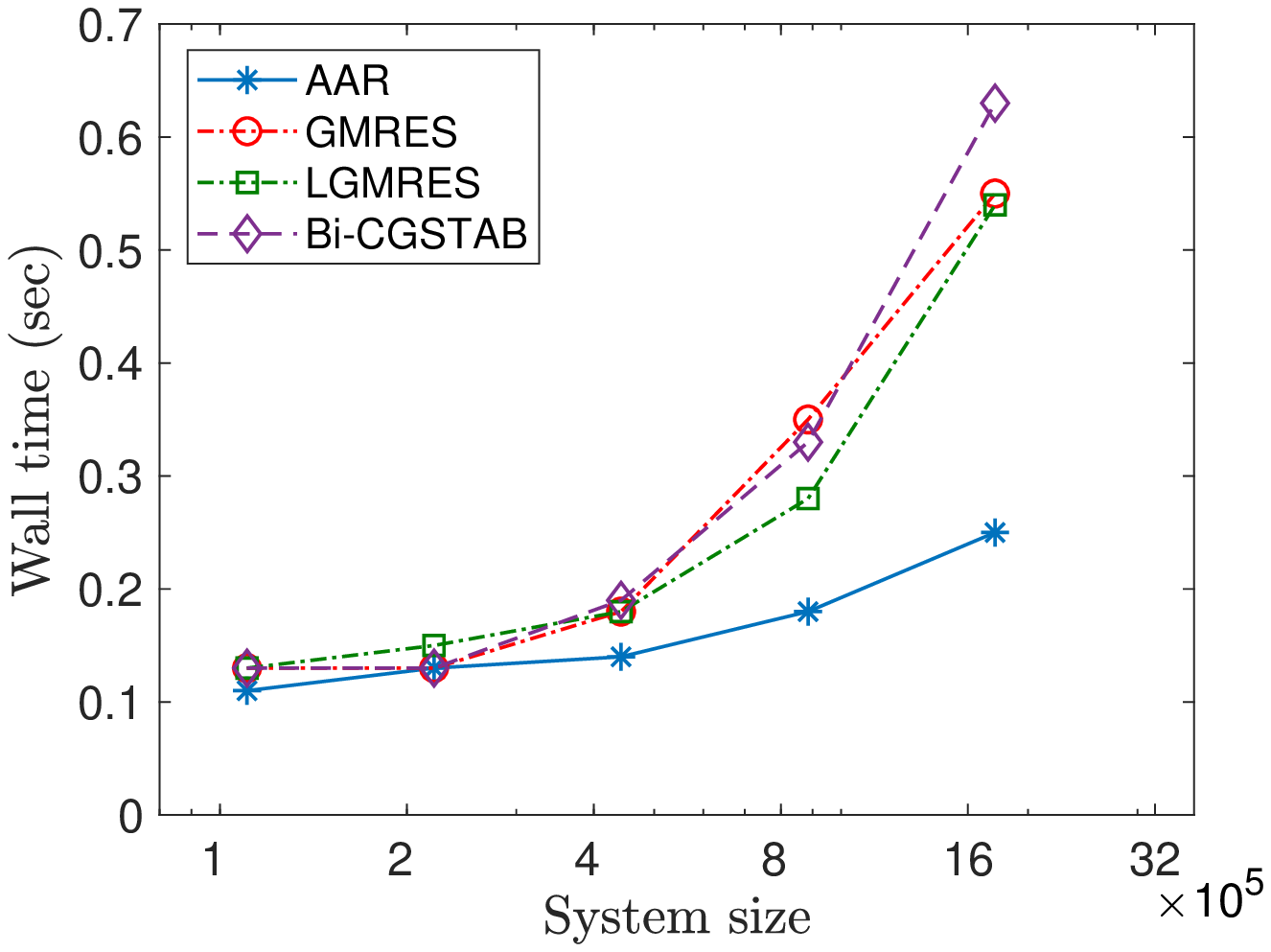}}   
\caption{Strong and weak scaling of AAR, GMRES, LGMRES, and Bi-CGSTAB for the periodic Helmholtz equation using block Jacobi preconditioning with ILU(0) factorization on each block in PETSc. In strong scaling, the minimum wall time taken by AAR is $1.38$, $1.43$, and $1.90$ times smaller than GMRES, LGMRES, and Bi-CGSTAB, respectively. In weak scaling, the CPU time taken by AAR, GMRES, LGMRES, and Bi-CGSTAB scales with system size as $\mathcal{O}(N^{1.28})$, $\mathcal{O}(N^{1.57})$, $\mathcal{O}(N^{1.51})$, $\mathcal{O}(N^{1.59})$, respectively.}
\label{Fig:Helmholtz}
\end{figure}

We first consider a $3 \times 3 \times 3$ aluminum supercell based on a face-centered cubic (FCC) unit cell having lattice constant $7.78$ Bohr, with atoms randomly displaced from ideal positions. We discretize the domain using a finite-difference grid with a mesh-size of $h=0.486$ Bohr, which is sufficient to achieve chemical accuracy in the energy and atomic forces. For the resulting linear system, we employ block Jacobi preconditioning with ILU(0) factorization on each block. Fig.~\ref{Fig:strong:h} shows the wall time taken by AAR, GMRES, LGMRES, and Bi-CGSTAB on 1, 8, 27, 64, 216, 512, and 1000 cores. We observe that even though the performances of all approaches are similar at low core counts (a consequence of requiring similar number of iterations), AAR starts demonstrating superior performance as the core count is increased. In particular, the minimum wall time achieved by AAR is a factor of $1.38$, $1.43$, and $1.90$ smaller than GMRES, LGMRES, and Bi-CGSTAB, respectively. This is a consequence of the significantly less global communication in AAR compared to the other methods.  

We next periodically replicate the above $3 \times 3 \times 3$ aluminum system along one direction by factors of $1$, $2$, $4$, $8$, and $16$, i.e., we generate $3 \times 3 \times 3$, $6 \times 3 \times 3$, $12 \times 3 \times 3$, $24 \times 3 \times 3$, and $48 \times 3 \times 3$ supercells. Correspondingly, we choose 64, 128, 256, 512, and 1024 computational cores. Again, we select $h = 0.486$ Bohr and employ block Jacobi preconditioning with ILU(0) factorization on each block. We plot the results of this weak scaling study in Fig.~\ref{Fig:weak:h}, from which we obtain $\mathcal{O}(N^{1.28})$, $\mathcal{O}(N^{1.57})$, $\mathcal{O}(N^{1.51})$, and $\mathcal{O}(N^{1.59})$ scaling with system size for AAR, GMRES, LGMRES, and Bi-CGSTAB, respectively. Notably, all approaches demonstrate slightly superlinear scaling even though the number of iterations remain constant. This is due to the increased cost of global communications at larger core counts. Therefore, the performance of AAR relative to the other methods is expected to further improve as the number of cores increases. 

\begin{figure}[h]
\centering
\subfloat[$h=0.216$ Bohr]{\label{Fig:app:mesh1}\includegraphics[width=0.47\textwidth]{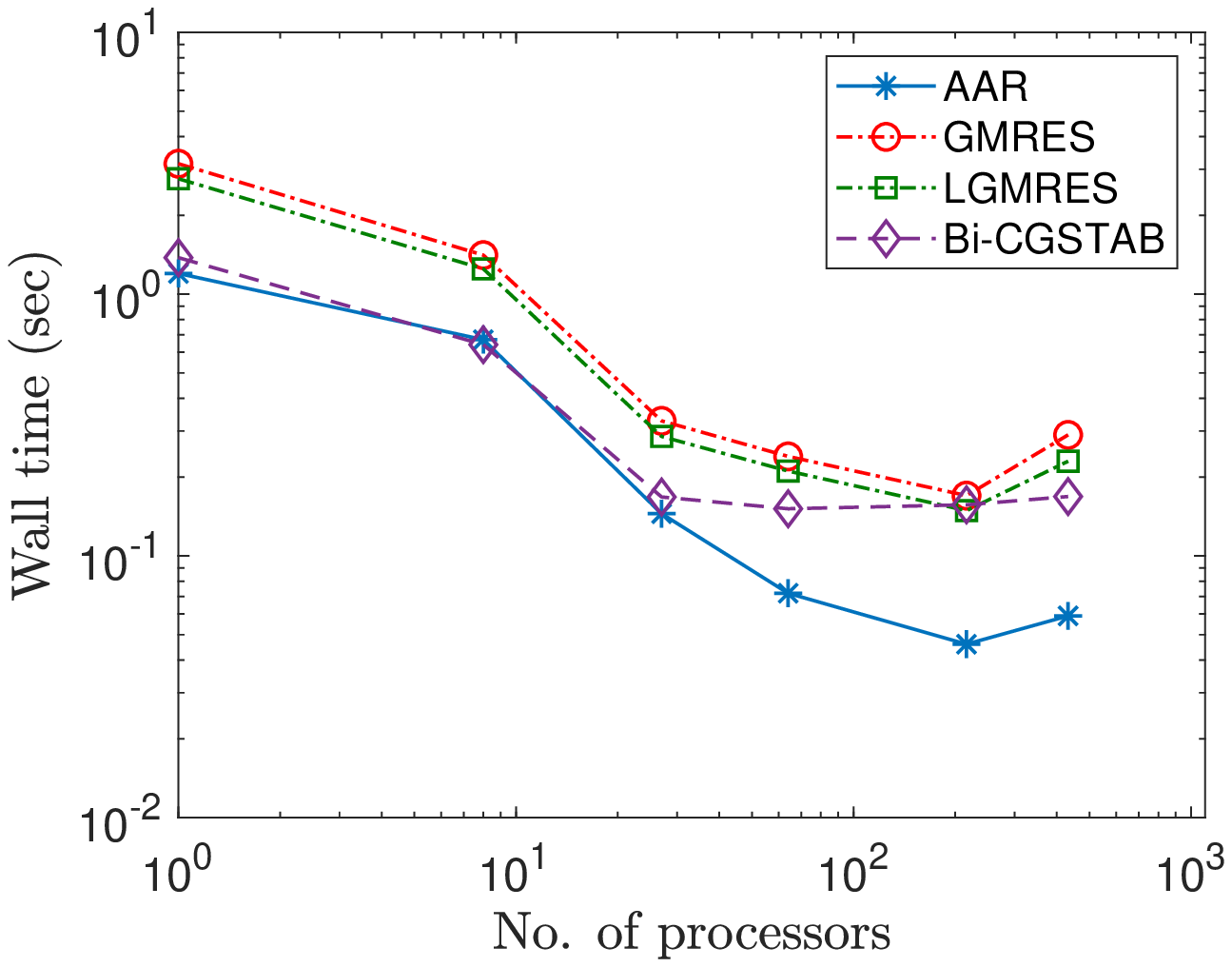}}
\subfloat[$h=0.108$ Bohr]{\label{Fig:app:mesh2}\includegraphics[width=0.47\textwidth]{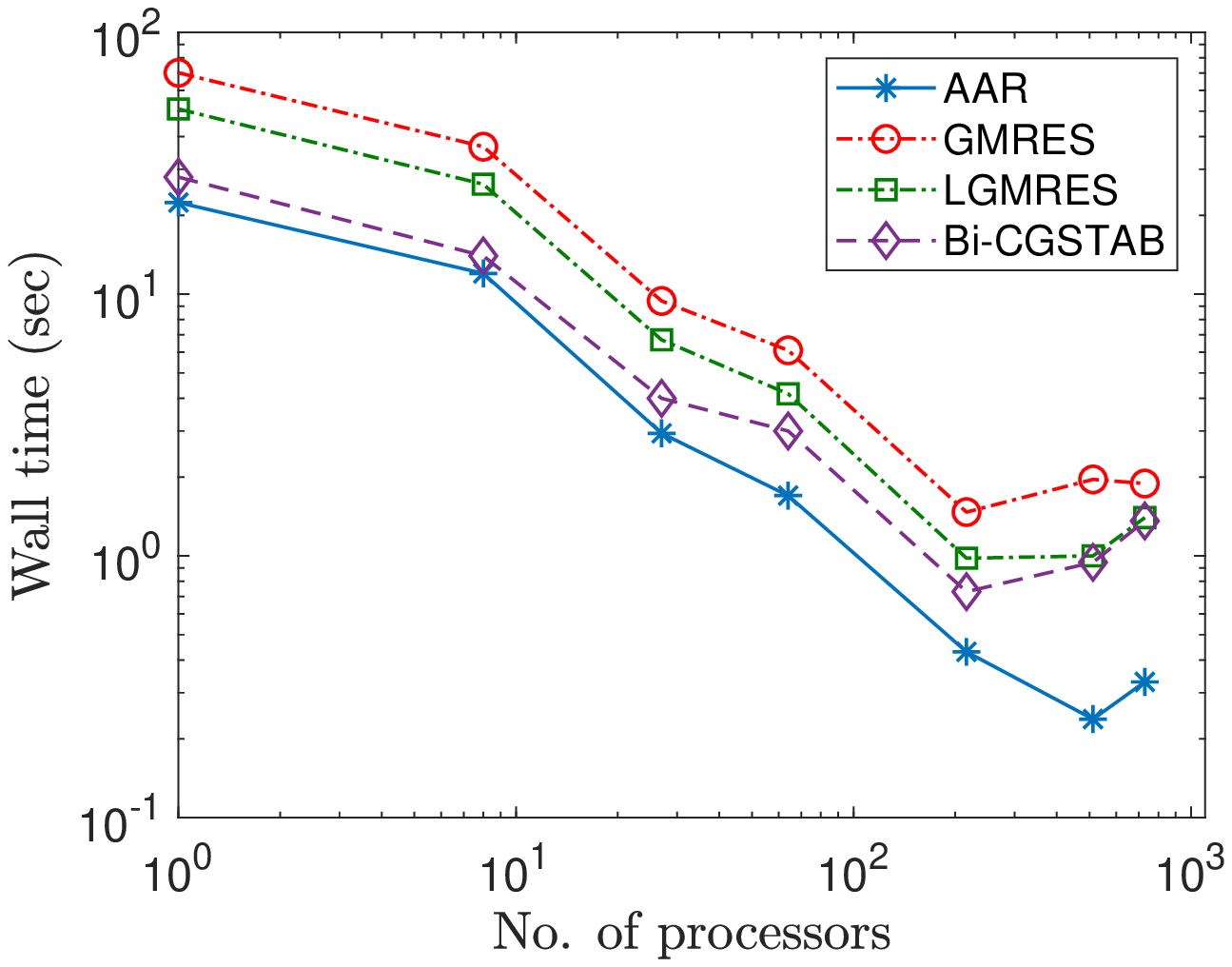}}
\caption{Strong scaling of AAR, GMRES, LGMRES, and Bi-CGSTAB for the periodic Helmholtz equation with Jacobi preconditioning in PETSc. For $h=0.216$ Bohr, the minimum wall time taken by AAR is $3.70$, $3.26$, and $3.28$ times smaller than GMRES, LGMRES, and Bi-CGSTAB, respectively. For $h=0.108$ Bohr, the minimum wall time taken by AAR is $6.13$, $4.08$, and $3.04$ times smaller than GMRES, LGMRES, and Bi-CGSTAB, respectively.}
\label{Fig:app:Helmholtz}
\end{figure}

It is worth noting that the performance gap between AAR and other approaches increases as the linear system becomes less well conditioned, as demonstrated in previous work for AAJ vs. GMRES in serial computations \cite{pratapa2016anderson}. In order to verify this result for AAR in parallel computations, we consider the Helmholtz equation for a $1 \times 1 \times 1$ Al supercell with lattice constant of $7.78$ Bohr and randomly displaced atoms. To demonstrate the effect of conditioning, we choose a simple Jacobi preconditioner $\mathbf{M}^{-1} = \mathbf{D}^{-1}$, where $\mathbf{D}$ is the diagonal part of $\mathbf{A}$. In Fig.~\ref{Fig:app:Helmholtz}, we present the strong scaling of AAR, GMRES, LGMRES, and Bi-CGSTAB for mesh-sizes of $h=0.216$ and $h=0.108$ Bohr. We observe that as the mesh gets finer and condition number of $\mathbf{A}$ becomes larger, the speedup of AAR over the other methods increases at both small and large core counts. Specifically, for $h=0.216$ Bohr, the minimum wall time taken by AAR is $3.70$, $3.26$, and $3.28$ times smaller than GMRES, LGMRES, and Bi-CGSTAB, respectively. The corresponding numbers for $h=0.108$ Bohr are $6.13$, $4.08$, and $3.04$, respectively. Therefore, we conclude that as the solution of the linear system becomes more challenging (e.g., in the absence of an effective preconditioner), the speedup of AAR over Krylov subspace approaches like GMRES and Bi-CGSTAB is expected to become more substantial in both the serial and parallel settings.\footnote{The upturns in strong scaling plots at largest core counts arise due to insufficient computational work per core relative to local inter-processor communications when the chosen computation is spread beyond a certain number of cores. This indicates the strong scaling limit of the current implementation for the chosen problem size.} 


\subsection{Density Functional Theory: Poisson equation} \label{sec:poisson}
Next, we study the relative performance of AAR and Conjugate Gradient (CG) methods for solving the periodic Poisson problem arising in real-space Density Functional Theory (DFT) calculations \cite{Pask2005,suryanarayana2013coarse, Pask2012,ghosh2016sparc2}:
\begin{equation} 
-\frac{1}{4\pi} \nabla^2 \phi(\br)  =  \rho(\br) + b(\br) \,\,\, \text{in} \,\,\, \Omega,\quad \begin{cases} V(\br)=V(\br+L_i\hat{\be}_i) \,\,\, \text{on} \,\,\, \partial\Omega \,, \\ \hat{\be}_i\cdot \nabla V(\br)=\hat{\be}_i\cdot \nabla V(\br+L_i\hat{\be}_i) \,\,\, \text{on} \,\,\, \partial\Omega \,, \end{cases}  \label{Eqn:poisson_p}  
\end{equation}
where $\phi(\br)$ is the electrostatic potential, $\rho(\br)$ is the electron density, $b(\br)$ is the pseudocharge density \cite{Pask2005,Phanish2010,Phanish2011,ghosh2016sparc}, and $\Omega$ is a cuboidal domain with side lengths $L_i$, unit vectors $\hat{\be}_i$ along each edge, and boundary $\partial \Omega$. The Poisson equation is discretized using sixth-order accurate finite-differences on a uniform grid with mesh-size $h=0.486$ Bohr, which is sufficient to achieve chemical accuracy in the energy and atomic forces. The computations are parallelized by decomposing the domain into cubical subdomains of equal size, with communication between processors handled via Message Passing Interface (MPI). 

In the Anderson extrapolation step of AAR, we again perform only one global communication call (i.e., \texttt{MPI\_Allreduce}) to simultaneously determine $r$ and the complete matrix $\mathbf{F}_k^{T}\mathbf{F}_k$, whose inverse in computed using the Moore-Penrose pseudoinverse. In all calculations, we again choose a vector of all zeros as the starting guess $\bx_0$, and a convergence tolerance of $\epsilon = 10^{-6}$ on the relative residual. Calculations on up to 1,024 cores were carried out on the same computer cluster as for the Helmholtz problem in Section~\ref{sec:helmholtz}. Larger calculations, up to 110,592 cores, were carried out on the Vulcan IBM BG/Q machine at the Lawrence Livermore National Laboratory, consisting of 24,576 compute nodes, with 16 computational cores and 16 GB memory per node, for a total of 393,216 cores and 1.6 PB memory.

\begin{figure}[h]
\centering
\subfloat[Strong scaling]{\label{Fig:strong:p}\includegraphics[width=0.47\textwidth]{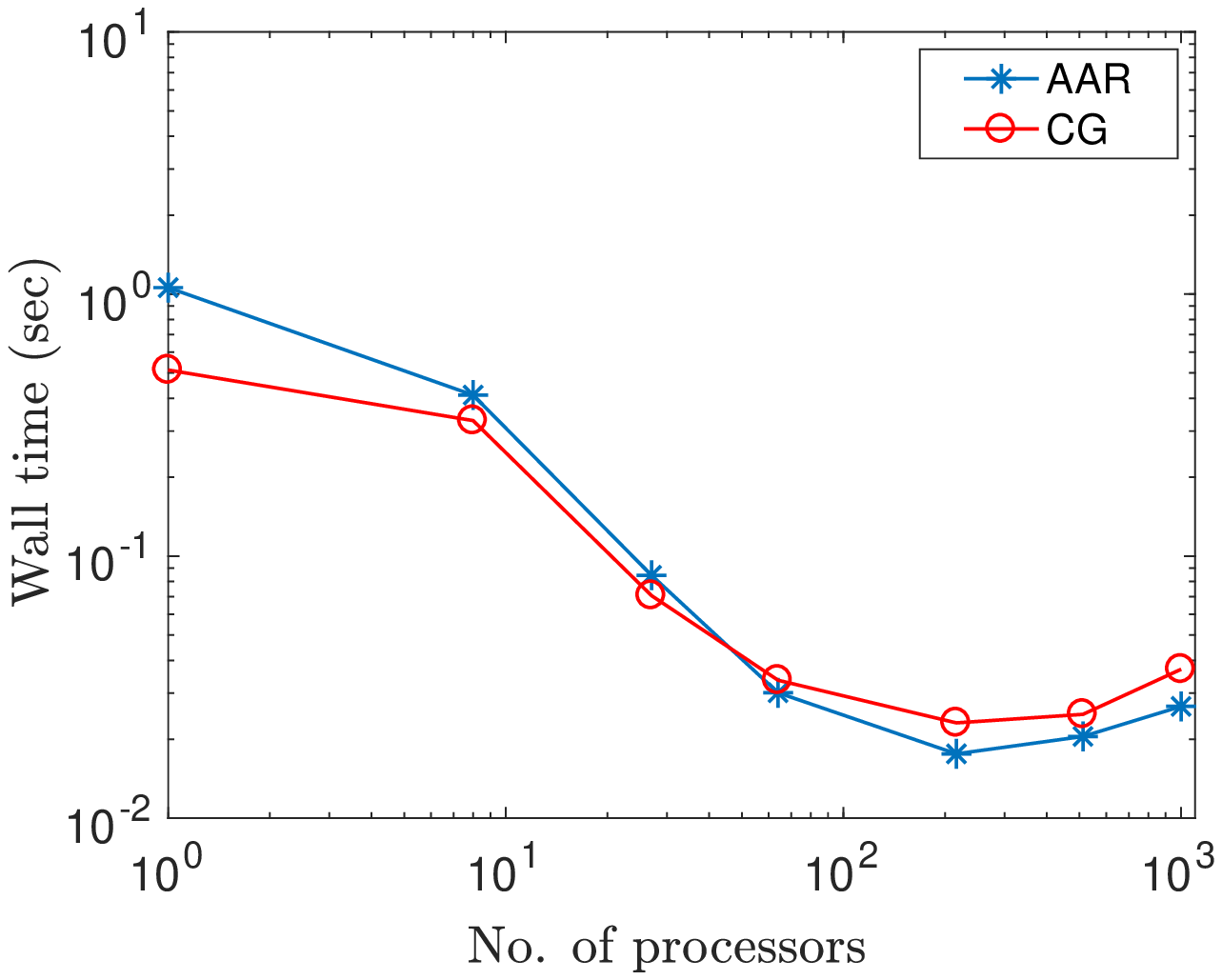}}
\subfloat[Weak scaling]{\label{Fig:weak:p}\includegraphics[width=0.47\textwidth]{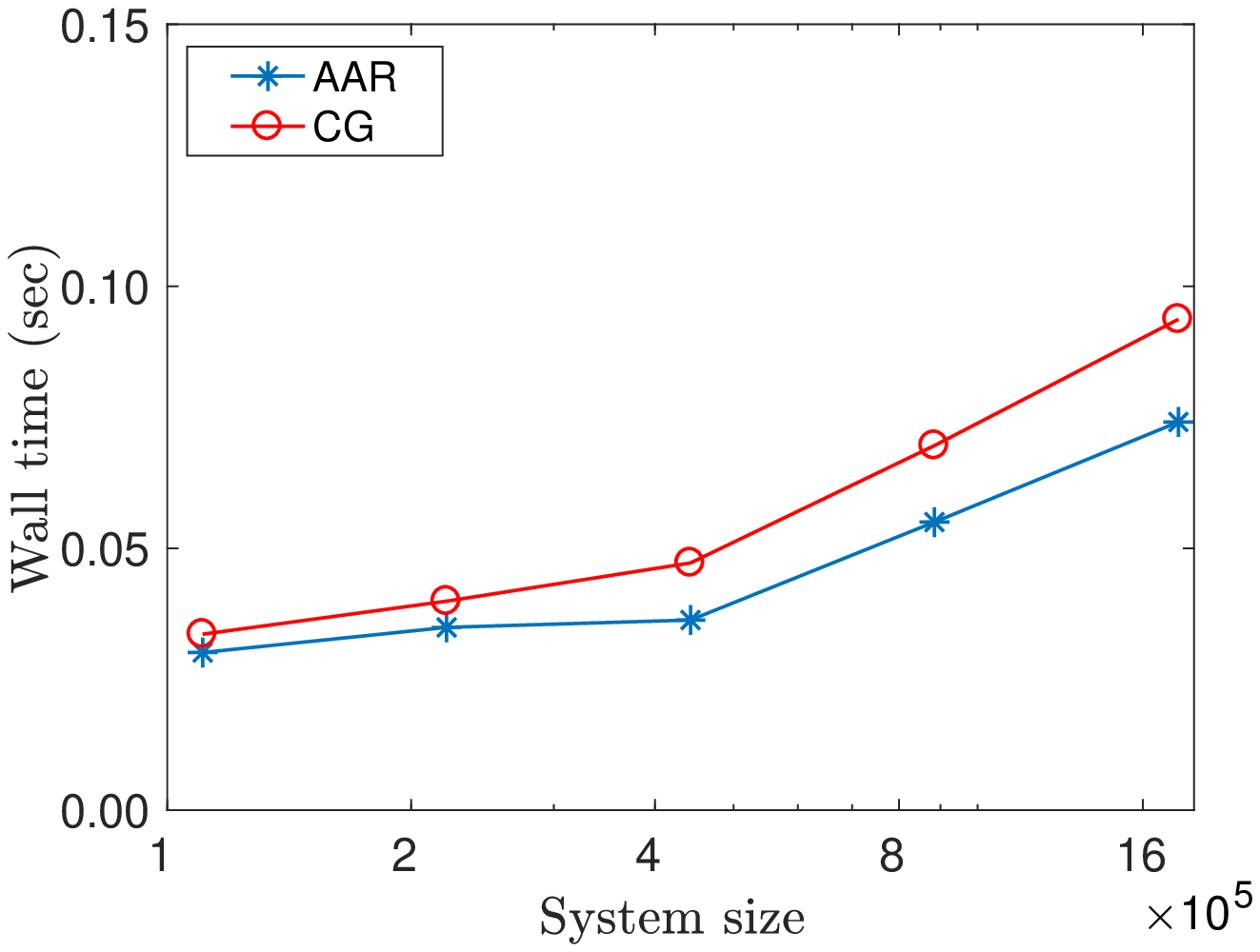}}
\caption{Strong and weak scaling of AAR and CG for the periodic Poisson equation using block Jacobi preconditioning with ILU(0) factorization on each block in PETSc. In strong scaling, the minimum wall time taken by AAR is $1.31$ times smaller than CG. In weak scaling, the CPU time taken by AAR and CG scales with system size as $\mathcal{O}(N^{1.33})$ and $\mathcal{O}(N^{1.38})$, respectively.}
\label{Fig:Poisson}
\end{figure}

We first consider a $3 \times 3 \times 3$ Al supercell based on a FCC unit cell with lattice constant 7.78 Bohr, with atoms randomly displaced from ideal positions, and again employ block Jacobi preconditioning with ILU(0) factorization on each block (default in PETSc) in the solution of the resulting linear systems. In Fig.~\ref{Fig:strong:p}, we plot the wall time taken by AAR and CG as implemented in PETSC\footnote{The PETSc implementation of AAR for real-valued systems is available as part of the code accompanying this paper. For CG, the inbuilt function in PETSc is utilized.} on 1, 8, 27, 64, 216, 512, and 1000 computational cores. At small core counts, CG demonstrates better performance than AAR by virtue of requiring fewer iterations to achieve convergence. However, as the number of cores is increased, the performance of AAR relative to CG improves, with the minimum wall time taken by AAR being a factor of $1.31$ smaller than CG by virtue of the lesser global communication required by AAR.

We next perform a weak scaling study by periodically replicating the above $3 \times 3 \times 3$ system along one direction by factors of $1$, $2$, $4$, $8$, and $16$, i.e., we generate $3 \times 3 \times 3$, $6 \times 3 \times 3$, $12 \times 3 \times 3$, $24 \times 3 \times 3$, and $48 \times 3 \times 3$ supercells. Correspondingly, we choose 64, 128, 256, 512, and 1024 computational cores. Again, we employ block Jacobi preconditioning with ILU(0) factorization on each block. In Fig.~\ref{Fig:weak:p}, we plot the wall time taken by AAR and CG for the resulting systems, from which we obtain the weak scaling with system size to be $\mathcal{O}(N^{1.33})$ and $\mathcal{O}(N^{1.38})$, respectively. As before, even though the number of iterations does not vary with system size, the increasing cost associated with global communications results in superlinear scaling for both approaches\footnote{Another factor contributing to the superlinear scaling for both AAR and CG is the inefficiency of communications between processors in the computer cluster used for the study.}. Therefore, the performance of AAR relative to CG is expected to further improve as core counts are increased, a result which we verify next. 

\begin{figure}[h]
\centering
\subfloat[Strong scaling]{\label{Fig:strong:p:llnl}\includegraphics[width=0.47\textwidth]{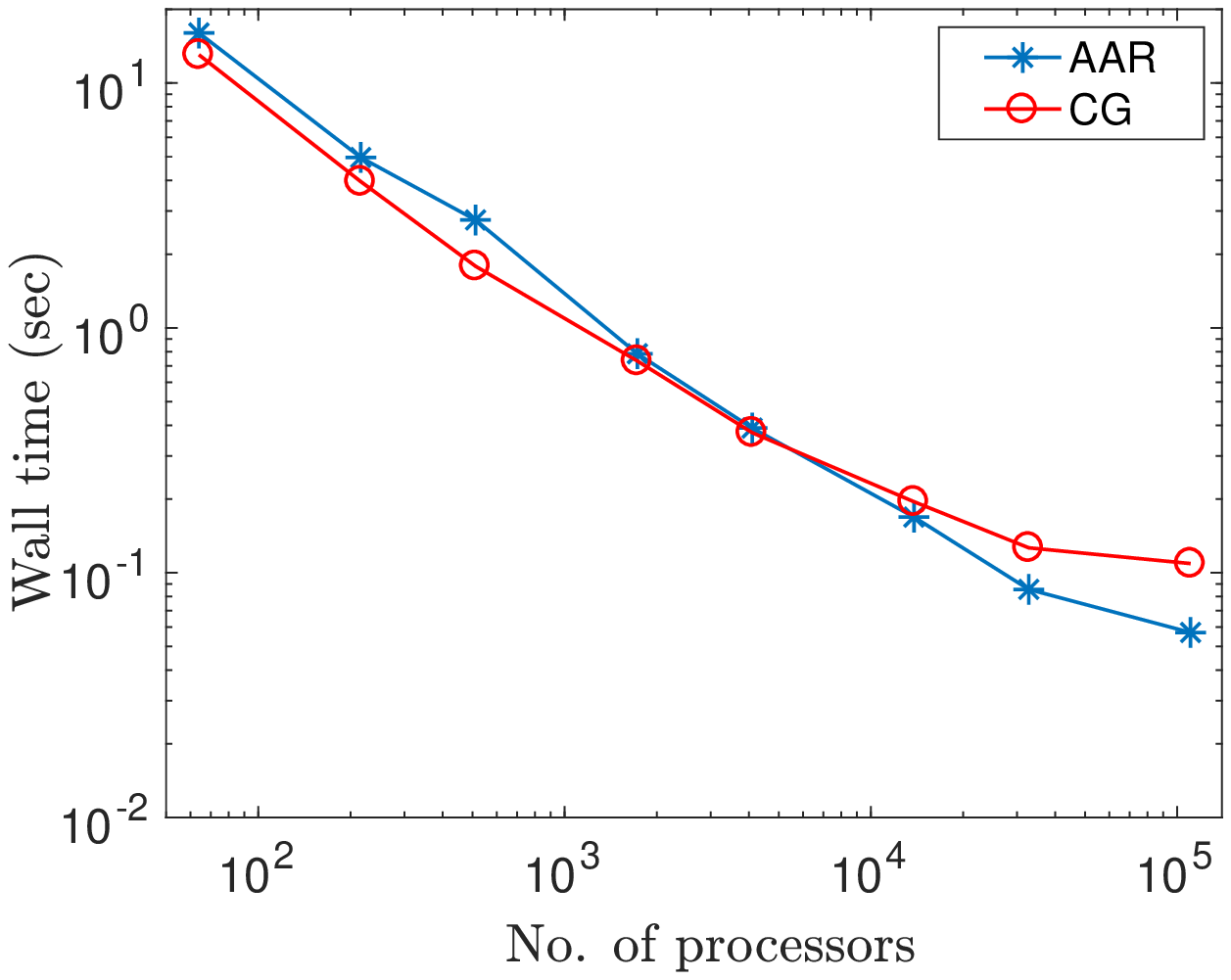}}
\subfloat[Weak scaling]{\label{Fig:weak:p:llnl}\includegraphics[width=0.47\textwidth]{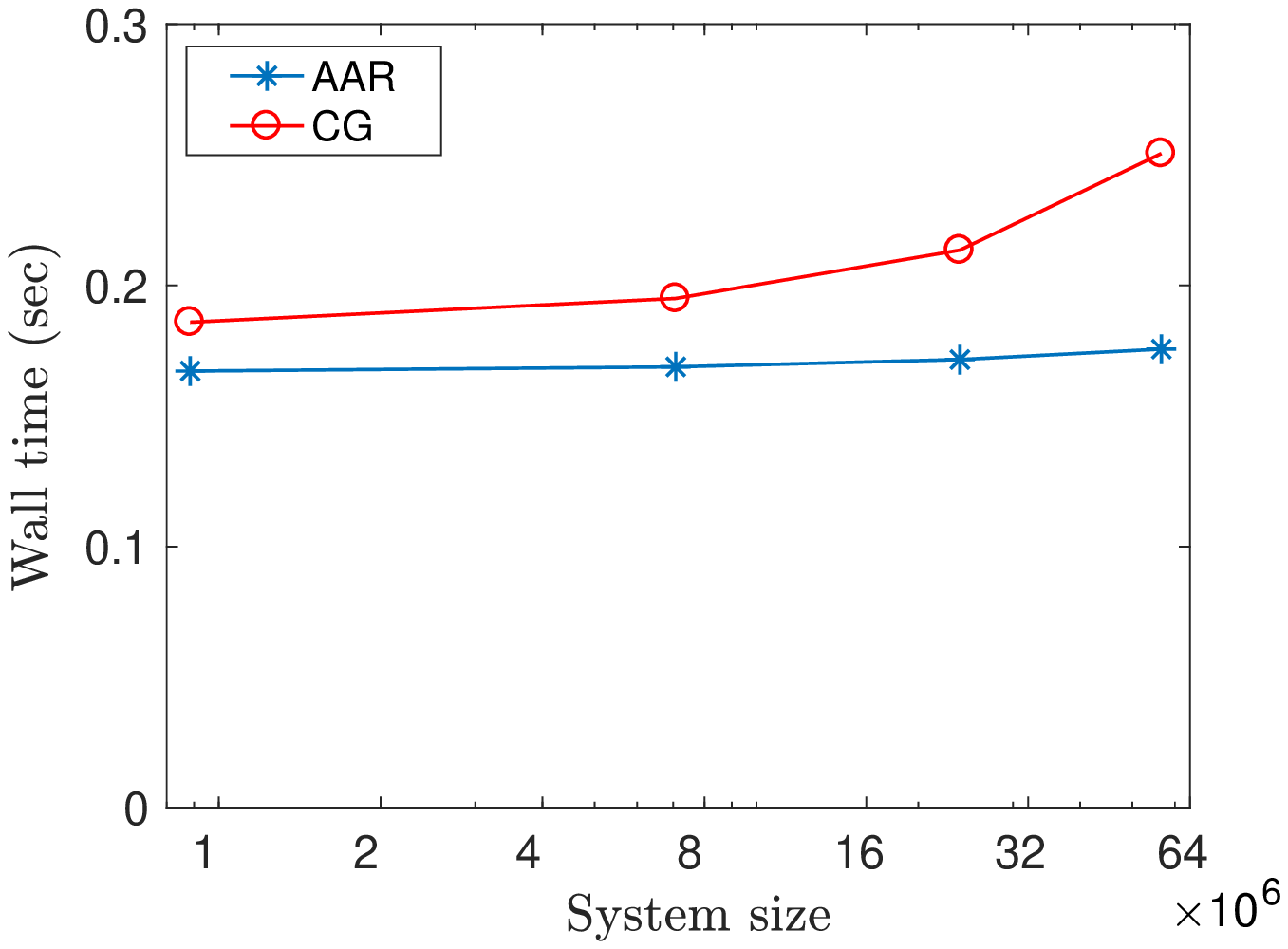}}
\caption{Strong and weak scaling of AAR and CG for the periodic Poisson problem with Jacobi preconditioning. In strong scaling, the minimum wall time taken by AAR is $1.91$ times smaller than CG. In weak scaling, the CPU time taken by AAR and CG scales with system size as $\mathcal{O}(N^{1.01})$ and $\mathcal{O}(N^{1.07})$, respectively.}
\label{Fig:Poisson:llnl}
\end{figure}

To assess the efficiency and scaling of AAR in larger-scale parallel calculations, up to 110,592 cores, we consider strong and weak scaling on the Vulcan IBM BG/Q machine at the Lawrence Livermore National Laboratory. We implement AAR and CG using C and MPI directly\footnote{The standalone implementation of AAR for real-valued systems using C and MPI directly is available as part of the code accompanying this paper.}, and choose a simple Jacobi preconditioner for parallel scalability. For the strong scaling study, we choose a $12\times 12 \times 12$ Al supercell with atoms randomly displaced, and a maximum of 110,592 computational cores. For the weak scaling study, we go from $6\times 6 \times 6$ supercell on 1728 processors to $24\times 24 \times 24$ supercell on 110,592 processors, with atomic displacements, electron density, and pseudocharge density periodically repeated from the $6\times 6 \times 6$ system. We present the results obtained in Fig.~\ref{Fig:Poisson:llnl}. We find that the minimum wall time achieved by AAR within the number of cores available for this study is $1.91$ times smaller than that achieved by CG. In addition, the weak scaling with system size for AAR and CG are $\mathcal{O}(N^{1.01})$ and $\mathcal{O}(N^{1.07})$, respectively. This demonstrates again the increased advantage of AAR over current state-of-the-art Krylov solvers in parallel computations as the number of processors is increased.

\section{Concluding remarks} \label{Sec:Conclusions}
We generalized the recently proposed Alternating Anderson-Jacobi (AAJ) method to include preconditioning and make it particularly well suited for scalable high-performance computing, and demonstrated its efficiency and scaling in the solution of large, sparse linear systems on parallel computers. Specifically, the AAR method employs Anderson extrapolation at periodic intervals within a preconditioned Richardson iteration to accelerate convergence while maintaining its underlying parallel scalability and simplicity to the maximum extent possible. 
In serial applications to nonsymmetric systems, we find that AAR is comparably robust to GMRES, using the same preconditioning, while often substantially outperforming it in time to solution; and find AAR to be more robust than Bi-CGSTAB for the problems considered. 
In parallel applications to the Helmholtz and Poisson equations, we find that AAR shows superior strong and weak scaling to GMRES, Bi-CGSTAB, and Conjugate Gradient (CG) methods, using the same preconditioning, with consistently shorter times to solution at larger processor counts. Finally, in massively parallel applications to the Poisson equation, on up to 110,592 processors, we find that AAR shows superior strong and weak scaling to CG, with shorter minimum time to solution.

Our findings suggest that the AAR method provides an efficient and scalable alternative to current state-of-the-art preconditioned Krylov solvers for the solution of large, sparse linear systems on high performance computing platforms, with increasing advantage as the number of processors is increased. Moreover, the method is simple and general, applying to symmetric and nonsymmetric systems, real and complex alike. Additional mathematical analysis which provides further insights into the performance of the AAR method and therefore enables the development of effective preconditioners tailored to it will enable still larger-scale applications, and so constitutes a potentially fruitful direction for future research.


\section{Acknowledgements}
This work was performed under the auspices of the U.S. Department of Energy by Lawrence Livermore National Laboratory under Contract DE-AC52-07NA27344. We gratefully acknowledge support from the Laboratory Directed Research and Development Program. P.S. and P.P. also gratefully acknowledge the support of National Science Foundation under Grant Number 1333500. 


\bibliographystyle{ReferenceStyle}
\bibliography{Manuscript}

\end{document}